# Reduction of abstract homomorphisms of lattices mod $p$ and rigidity

Chandrashekhar Khare & Dipendra Prasad

## 1 Introduction

Let $A$ be an abelian variety over a number field $K$. By the Mordell-Weil theorem the group $A(K)$ of $K$-rational points of $A$ is a finitely generated abelian group. There is a finite set $S$ of places of $K$ such that there is an abelian scheme $\mathcal{A}$ over $\mathcal{O}_S$, the ring of $S$-integers of $K$, with generic fibre $A$. For $v$ not in $S$ (all the places $v$ we consider below are *outside* $S$ even if this is not mentioned explicitly) we denote by $A_v$ the fibre of $\mathcal{A}$ at $v$, $k_v$ the residue field at $v$. Consider the specialisation map $sp_v : A(K) \to A_v(k_v)$ and denote by $A(v)$ its image. We say that a homomorphism $\phi : A(K) \to A(K)$ specialises mod $v$ if there is a homomorphism $\phi_v : A(v) \to A(v)$ such that the diagram

$$\begin{array}{ccc} A(K) & \xrightarrow{\phi} & A(K) \\ sp_v \downarrow & & \downarrow sp_v \\ A(v) & \xrightarrow{\phi_v} & A(v) \end{array}$$

commutes. The following question is the starting point of this paper.

**Question 1** *Let*

$$\phi : A(K) \to A(K)$$

*be a homomorphism of abelian groups that specialises mod $v$ for almost all place $v$ of $K$ to give a homomorphism $\phi_v$ as above. Is the restriction of $\phi$ to a subgroup of finite index of $A(K)$ induced by an endomorphism $\alpha_\phi \in \mathrm{End}_K(A)$?*

This question is inspired by the results of the paper [RS] and the study of compatible systems of mod $p$ Galois representations in [K]. Indeed the hypotheses of the question can be viewed as the giving of a compatible system $\{\phi_v\}$ of homomorphisms of $A(v)$ and then one would like to prove a *reciprocity law* saying that these arise in the only natural way, i.e., from endomorphisms of $A$. This can also be considered as a *rigidity* theorem in analogy with similar conclusions for homomorphism of lattices (i.e. discrete subgroups of Lie groups with finite volume for the quotient) in semi-simple Lie groups. The following theorem responds to the question when the only endomorphisms of $A$ (over $\overline{K}$) are multiplication by $n$ and when dimension of $A$ is either 2, 6 or odd.

**Theorem 1** *Let $A$ be an abelian variety over a number field $K$ and assume that $\mathrm{End}(A) = \mathbf{Z}$ and $\dim(A) = g$ is either odd or $g = 2$ or $6$.*

1. *Then for any abstract homomorphism*

$$\phi : A(K) \to A(K)$$

   *that specialises mod $v$ for almost all places $v$ of $K$, the restriction of $\phi$ to a subgroup of $A(K)$ of finite index is given by the multiplication-by-n map for some integer $n$.*

2. *In fact, if $A(K)$ is not torsion then $\phi$ itself is induced by the multiplication-by-n map for some integer $n$.*

The restrictions on dimension arise from those of Serre's theorems in article 137 of [S-IV]. We prove theorem 1 in Section 2.1. In Section 2.2 we prove a general lemma that together with results of [RS] allows us to answer Question 1 in the case of CM elliptic curves. It is concievable that the methods used in the proof of the theorem might be extended to the case when $A$ is absolutely simple but is allowed to have non-trivial endomorphisms. But the restriction that $A$ is (absolutely) simple is a serious restriction, and when $A$ is not simple our methods are found wanting. Thus we are unable to answer Question 1 in the simplest case when $A$ is not simple, e.g. in the case $A = E_1 \times E_2$ where $E_1$ and $E_2$ are elliptic curves over $K$. To be specific, we are unable to answer the following very particular case of Question 1 implicit in whose affirmative answer is a novel criterion for elliptic curves to be isogenous. We thank Bas Edixhoven for this formulation.

**Question 2** *Let $E_1, E_2$ be elliptic curves over $K$ and consider a homomorphism $\phi : E_1(K) \to E_2(K)$ that specialises for almost all places $v$ of $K$ to give a homomorphism $\phi_v : E_1(v) \to E_2(v)$. Then is the homomorphism $\phi$, restricted to a subgroup of finite index of $E_1(K)$, induced by an isogeny from $E_1$ to $E_2$? One can in fact ask the following stronger question. Let $P$ (resp. $Q$) be a non-torsion point in $E_1(K)$ (resp. $E_2(K)$) such that for almost all places $v$ of $K$, the order of $Q$ mod $v$ divides the order of $P$ mod $v$ then are $E_1$ and $E_2$ are isogenous, and $Q$ related to $P$ by an isogeny from $E_1$ to $E_2$?*

We discuss this question in Section 3 by drawing a line between it and the conjectures of Lang-Trotter about analogs for elliptic curves of Artin's conjecture on primitive roots. In Section 4 we answer an analog of Question 1 for subgroups of finite index of $SL_2(\mathbf{Z})$, and in Section 5 consider analogs for rational varieties.

## 2 Rigidity for abelian varieties

### 2.1 Proof of Theorem 1

For the sake of exposition we prefer to break up the proof of the theorem into 2 parts which correspond to the 2 parts in the statement of it. The proof is directly inspired and follows very closely the methods of [RS]. The novelty of the proof is only in the proofs of Lemma 3 and 5, that replace the use of Siegel's theorem on finiteness of integral points on elliptic curves to prove their analogs in [RS].

The following lemma is well known, cf. S.Lang's book, Algebra, section 10 of the chapter on Galois theory for $i = 1$. It also follows from generalities on cohomology once we know it is true for $i = 0$, which is of course clear.

**Lemma 1** *Let $G$ be a group, and $E$ a $G$-module. Let $\tau$ be an element in the center of $G$. Then $H^i(G, E), i = 0, 1, \ldots$ is annihilated by the map $x \to \tau x - x$.*

**Lemma 2** *Let $G_{\ell^n} = \mathrm{Gal}(K_{\ell^n}/K)$ with $K_{\ell^n} = K(A[\ell^n])$, and $G_{\ell^\infty} = \mathrm{Gal}(K_{\ell^\infty}/K)$ with $K_{\ell^\infty} = \cup_n K_{\ell^n}$. Then $H^1(G_{\ell^m}, A[\ell^n])$ $(m \geq n)$ and $H^1(G_{\ell^\infty}, A[\ell^n])$ are of finite orders, bounded independent of $m$ and $n$.*

**Proof :** We note that $G_{\ell^\infty}$ being a compact $\ell$-adic Lie group, is topologically finitely generated, hence each finite quotient such as $G_{\ell^m}$ is generated by a set of elements of cardinality independent of $m$.

From the definition of $H^1(G, A)$ in terms of maps $\phi$ from $G$ to $A$ such that $\phi(g_1 g_2) = \phi(g_1) + g_1 \phi(g_2)$, it follows that an element of $H^1(G, A)$ is determined by a map on a set of generators of $G$. Since $A[\ell^n] \cong (\mathbf{Z}/\ell^n)^{2g}$ as abelian groups, it follows that $H^1(G_{\ell^m}, A[\ell^n])$ is a finitely generated abelian group which is generated by a set of elements of cardinality independent of $m, n$.

It follows from a theorem of Bogomolov that the $\ell$-adic Lie group $G_{\ell^\infty}$ contains homotheties congruent to 1 modulo $\ell^N$ for some integer $N > 0$. Therefore by Lemma 1, $H^1(G_{\ell^m}, A[\ell^n])$ is annihilated by $\ell^N$. It follows that $H^1(G_{\ell^m}, A[\ell^n])$ is a finitely generated abelian group which is generated by a set of elements of cardinality independent of $m, n$, and annihilated by $\ell^N$, and thus is of finite order, bounded independent of $m, n$.

The statement about $H^1(G_{\ell^\infty}, A[\ell^n])$ follows either by noting that the cohomology $H^1(G_{\ell^\infty}, A[\ell^n])$ can be calculated in terms of continuous cochains on $G_{\ell^\infty}$, for which the earlier argument applies as well, or by noting that $H^1(G_{\ell^\infty}, A[\ell^n])$ is the direct limit of $H^1(G_{\ell^m}, A[\ell^n])$ (direct limit over $m$), and direct limit of finitely generated abelian groups each of which is generated by a set of elements of cardinality independent of $n$, and each annihilated by $\ell^N$, is of order bounded independent of $n$.

**Lemma 3** *Given an abelian variety $A$ over $K$ and a point $P \in A(K)$ of infinite order and any prime $\ell$, there are infinitely many $v$'s (in fact a positive density of such) such that the reduction of $P$ mod $v$ has order divisible by $\ell$.*

**Proof:** We claim that there is a sufficiently large $n$ such that the extension $K_{P,\ell^n} = K(A[\ell^n], \frac{1}{\ell^n}.P)$ is a non-trivial extension of $K_{\ell^n} = K(A[\ell^n])$. If we grant the claim then there is a positive density of $v$ that split in $K(A[\ell^n])$ (we denote the Galois group $\mathrm{Gal}(K_{\ell^n}/K)$ by $G_{\ell^n}$) but not in $K(A[\ell^n], \frac{1}{\ell^n}.P)$, and by inspection we see that for such $v$'s the reduction of $P$ mod $v$ has order divisible by $\ell$.

From Lemma 2, it follows that $H^1(G_{\ell^n}, A[\ell^n])$ is bounded independently of $n$. Using maps between the Kummer sequences

$$
\begin{CD}
0 @>>> A(K)/\ell^n A(K) @>>> H^1(G_K, A[\ell^n]) @>>> H^1(G_K, A)[\ell^n] @>>> 0 \\
@. @VVV @VVV @VVV \\
0 @>>> A(K_{\ell^n})/\ell^n A(K_{\ell^n}) @>>> H^1(G_{K_{\ell^n}}, A[\ell^n]) @>>> H^1(G_{K_{\ell^n}}, A)[\ell^n] @>>> 0
\end{CD}
\tag{1}
$$

the claim follows for large enough $n$, putting together the observation that the kernel of the restriction map $H^1(G_K, A[\ell^n]) \to H^1(G_{K_{\ell^n}}, A[\ell^n])$ is $H^1(G_{\ell^n}, A[\ell^n])$, and the fact that as $P$ is non-torsion the image of $P$ under the coboundary map (in the first exact sequence) in $H^1(G_K, A[\ell^n])$ has unbounded order as $n$ varies. From this the claim follows and hence the lemma.

**Lemma 4** *(Lemma 4.2 of [RS]) Let $G$ be a subgroup of $GL_n(\mathbf{F}_\ell)$ and let $H_1$ and $H_2$ be 2 copies of $\mathbf{F}_\ell^n$ with the natural action of $G$. Consider $\mathbf{G}$ which is the semidirect product of $H_1 \oplus H_2$ by $G$. Let $\sigma$ be an element in $G$, such that for any $h_1 \in H_1$, $(h_1, \sigma)$ is conjugate to $(h_2, \tau)$ in $\mathbf{G}$ for some $h_2 \in H_2, \tau \in G$. (We denote an element $(h_1, 0)$ or $(0, h_2)$ of $H_1 \oplus H_2$ by $h_1$ or $h_2$.) Then $\sigma$ cannot have 1 as eigenvalue.*

**Proof:** The proof is identical to the $n = 2$ case done in [RS], and we reproduce it here for the reader's convenience. A simple calculation yields that
$$(h', \rho)(h_1, \sigma)(h', \rho)^{-1} = (\rho^{-1}(h_1 + (\sigma - 1)h'), \rho\sigma\rho^{-1})$$
for $(h', \rho) \in \mathbf{G}$ and from this we deduce that for every $h_1 \in H_1$ there exists $h' \in H_1 \oplus H_2$ such that $h_1 + (\sigma - 1)h' \in H_2$. This implies $H_1 \subset (\sigma - 1)(H_1 \oplus H_2) + H_2$, thus $H_1 \oplus H_2 \subset (\sigma - 1)H_1 \oplus H_2$ and hence $\sigma - 1$ is surjective in its action on $H_1$.

We have the following key proposition whose proof is exactly as in Section 3 of [RS].

**Proposition 1** *If $P$ and $Q$ are non-torsion points $\in A(K)$ such that for almost all $v$ the order of $Q$ mod $v$ (i.e., the order of its image in $A_v(k_v)$) divides the order of $P$ mod $v$ (i.e., the order of its image in $A_v(k_v)$), then $aP = bQ$ for some $a, b \in \mathbf{Z}$.*

**Proof:** We consider only primes $\ell$ so that denoting the image of the action of $G_K$ on $A[\ell]$ by $G_\ell$, $H^1(G_\ell, A[\ell])$ is trivial, and further such that $G_\ell =$ $GSp_{2g}(\mathbf{F}_\ell)$. Under our assumptions and by the theorems of Serre in [S-IV] this is true for all $\ell$ large enough (see Théoréme C on page 40 and the Corollaire on page 51 of the letter to Vigneras in [S-IV]). We further ensure at the cost of excluding a further finite set of primes $\ell$, and bearing in mind that $P$ and $Q$ are assumed to be of infinite order, that $P$ and $Q$ are not divisible by $\ell$ in $A(K)$.

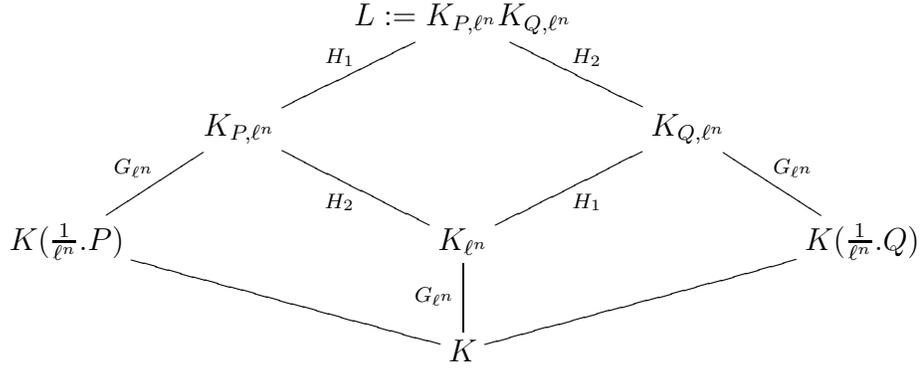

(Consider $n = 1$ in the above picture for the ensuing arguments.) We consider the fields $K_{P,\ell}$ and $K_{Q,\ell}$ over $K$, and want to prove that $K_{P,\ell} = K_{Q,\ell}$. If not (as we assume only to derive a contradiction later) by the irreducibility of $A[\ell]$ as a $G_K$ module we see that these fields are linearly disjoint over $K_\ell$ and their Galois groups over $K_\ell$ are both isomorphic to $A[\ell] \simeq \mathbf{F}_\ell^{2g}$. We consider $R$ (resp., $R'$) such that $\ell.R = P$ (resp., $\ell.R' = Q$) and the field $K(R)$ (resp., $K(R')$): note that by our assumptions on $\ell$ it follows that $\mathrm{Gal}(K_{P,\ell}/K(R)) \simeq G_\ell$ (resp., $\mathrm{Gal}(K_{Q,\ell}/K(R')) \simeq G_\ell$). Thus $\mathrm{Gal}(K_{P,\ell}/K)$ (resp., $\mathrm{Gal}(K_{Q,\ell}/K)$) is the semidirect product of $H_2 := \mathrm{Gal}(K_{P,\ell}/K_\ell)$ (resp., $H_1 := \mathrm{Gal}(K_{Q,\ell}/K_\ell)$, both isomorphic to $A[\ell]$) by $\mathrm{Gal}(K_{P,\ell}/K(R))$ (resp., $\mathrm{Gal}(K_{P,\ell}/K(R'))$, both isomorphic to $G_\ell$). Let $L = K_{P,\ell}K_{Q,\ell}$ and $\Omega := \mathrm{Gal}(L/K)$ which is the semidirect product (of the direct product) $H_1 \oplus H_2$ by $G_\ell$.

Consider an element $\sigma$ in $G_\ell$ that has 1 as an eigenvalue with multiplicity one and which exists by our assumptions on $\ell$. Then we claim that for any $h_1 \in H_1$, $(h_1, \sigma) \in H_1 \times G_\ell \simeq \mathrm{Gal}(L/K(R))$ is conjugate in $\Omega$ to an element $(h_2, \tau) \in H_2 \times G_\ell$. This will give the contradiction that $\sigma - 1$ is invertible by Lemma 4 thus proving our assumption wrong and hence that $K_{P,\ell} = K_{Q,\ell}$. To see the claim consider a degree one prime $v$ of $K$ whose Frobenius conjugacy class in $\Omega$ is the conjugacy class of $(h_1, \sigma)$ (what we have to show is that this conjugacy class intersects $H_2 \times G$). Thus the $\ell$-part of $A_v(k_v)$ is non-trivial and cyclic. Further, $v$ has a degree one prime above it in $K(R)$, and thus the order of $P$ in $A_v(k_v)$ divides $m/\ell$ where $m = |A_v(k_v)|$. Thus the order of $Q$ in $A_v(k_v)$ also divides $m/\ell$ by the hypothesis of the proposition. But this by our assumption on $\sigma$, that it has the eigenvalue 1 with multiplicity 1 (which as already noted implies that the $\ell$-part of $A_v(k_v)$ is cyclic and non-trivial), yields that $Q$ is divisible by $\ell$ in $A_v(k_v)$ and thus by the transitivity of the action of $G_\ell$ on $A[\ell]$, $v$ has a degree one prime above it in $K(R')$ which proves our claim, i.e., the Frobenius conjugacy class of $v$ contains $(h_1, \sigma)$ and intersects the semidirect product $H_2 \times G \simeq \mathrm{Gal}(L/K(R'))$. Thus we have proved that for almost all primes $\ell$, $K_{P,\ell} = K_{Q,\ell}$ and hence by Theorem 3.2 in [JR] we get that $aP = bQ$ for some integers $a, b$.

**Proof of part 1 of Theorem 1:** We claim that $b$ divides $a$. If not there would be a prime $\ell$ which divides $b$ to a higher power than the highest power (say $\ell^n$) to which it divides $a$. Then use Lemma 3 and consider the infinitely many places $v$ such that the order of the image of $\ell^n.P \in A_v(k_v)$ is divisible by $\ell$. This gives a contradiction to our hypothesis that the order of $Q$ divides the order of $P$ in $A_v(k_v)$ for almost all $v$. Thus $\phi(P) = Q$ is a multiple of $P$ up to a torsion point. Let $t$ be an integer that annihilates the torsion of $A(K)$ and consider the free abelian group $tA(K)$. For each $P \in tA(K)$ we see that $\phi(P)$ is not torsion, and thus that that $\phi(P) = mP + R$ for some $m$ that might a priori depend on $P$ and for $R$ a torsion point. But then $R$ is zero as both $\phi(P)$ and $P$ are in the torsion-free subgroup $tA(K)$. Now considering $P + Q$ for $P, Q \in tA(K)$ it follows that $m$ is independent of $P$.

**Proof of part 2 of Theorem 1:** To complete the proof of part 2 of Theorem 2 we need the following lemma.

**Lemma 5** *Given an abelian variety $A_{/K}$ that satisfied hypotheses of Theorem 1 and a point $P \in A(K)$ of infinite order and any prime $\ell$, there are infinitely many $v$'s (in fact a positive density of such) such that the reduction of $P$ mod $v$ has order prime to $\ell$.*

**Proof:** From the proof of Lemma 3, it follows that $K_{P,\ell^\infty}$ is an infinite extension of $K_{\ell^\infty}$. Under the hypothesis of this lemma, $G_{\ell^\infty}$ is an open subgroup of finite index of $GSp_{2g}(\mathbf{Z}_\ell)$ by a theorem of Serre [S-IV]. (If the abelian variety is not principally polarised, the image of the Galois group does not necessarily sit inside a symplectic similitude group, and hence a slight modification to the argument given below will have to be done which we leave to the reader.) If we denote by $E_{\ell^\infty}$ the Galois group of $K_{P,\ell^\infty}$ over $K$, and $A_{\ell^\infty}$, the Galois group of $K_{P,\ell^\infty}$ over $K_{\ell^\infty}$, we have an exact sequence of groups,

$$0 \to A_{\ell^\infty} \to E_{\ell^\infty} \to G_{\ell^\infty} \to 1.$$

This exact sequence is a subsequence of the following (split) exact sequence of $\mathbf{Z}_\ell$-rational points of algebraic groups:

$$0 \to \mathbf{Z}_\ell^{2g} \to E(\mathbf{Z}_\ell) \to GSp(2g, \mathbf{Z}_\ell) \to 1,$$

(i.e., $A_{\ell^\infty} \subset \mathbf{Z}_\ell^{2g}, E_{\ell^\infty} \subset E(\mathbf{Z}_\ell), G_{\ell^\infty} \subset GSp(2g, \mathbf{Z}_\ell)$), in which $E(\mathbf{Z}_\ell)$ is the semi-direct product $GSp(2g, \mathbf{Z}_\ell)$ with $\mathbf{Z}_\ell^{2g}$. The embedding of $E_{\ell^\infty}$ inside $E(\mathbf{Z}_\ell)$ is obtained by choosing a sequence of points $P_n$ in $A$ with $\ell \cdot P_1 = P$, and $\ell \cdot P_{i+1} = P_i$. Since the action of $GSp(2g, \mathbf{Z}_\ell)$ on $\mathbf{Z}_\ell^{2g}$ is irreducible, and $G_{\ell^\infty}$ is an open subgroup of $GSp(2g, \mathbf{Z}_\ell)$, it follows that $A_{\ell^\infty}$ is an open subgroup of $\mathbf{Z}_\ell^{2g}$, and hence $E_{\ell^\infty}$ is an open subgroup of the semi-direct product $GSp(2g, \mathbf{Z}_\ell)$ with $\mathbf{Z}_\ell^{2g}$.

We will prove that the intersection of the fields $K_{P,\ell^n}$ and $K_{\ell^{n+1}}$ is $K_{\ell^n}$ for some $n$ large enough. This, together with the theorem of Bogomolov recalled earlier, will imply that there is a positive density of primes in $K$ which are split in $K_{P,\ell^n}$ and for which the Frobenius as an element of $\mathrm{Gal}(K_{\ell^{n+1}}/K)$ is a non-trivial homothety in $GSp(2g, \mathbf{Z}/\ell^{n+1})$ which is congruent to 1 modulo $\ell^n$. For such primes $v$, it can be easily seen that the order of $P$ modulo $v$ is not divisible by $\ell$, completing the proof of the lemma.

It thus suffices to prove that the intersection of the fields $K_{P,\ell^n}$ and $K_{\ell^{n+1}}$ is $K_{\ell^n}$.

Let $E_{\ell^n}$ (resp. $G_{\ell^n}$) denote the Galois group of $K_{P,\ell^n}$ (resp. $K_{\ell^n}$) over $K$, and let $A_{\ell^n}$ denote the Galois group of $K_{P,\ell^n}$ over $K_{\ell^n}$. We have the exact sequence of groups,

$$0 \to A_{\ell^n} \to E_{\ell^n} \to G_{\ell^n} \to 1.$$

It is clear that the intersection of the fields $K_{P,\ell^n}$ and $K_{\ell^{n+1}}$ is $K_{\ell^n}$ if and only if inside the group $E_{\ell^n}$ which is a quotient of $E_{\ell^{n+1}}$, the image of $A_{\ell^{n+1}}$ is $A_{\ell^n}$.

Since $E_{\ell^\infty}$ is an open subgroup of the semi-direct product $GSp(2g, \mathbf{Z}_\ell)$ with $\mathbf{Z}_\ell^{2g}$, it follows that $E_{\ell^\infty}$ contains the natural congruence subgroup of level $\ell^n$ in this semi-direct product for some $n$, say $n = n_o$. From this it is clear that $E_{\ell^{n+1}}$ is the full inverse image of $E_{\ell^n}$ under the natural mapping from $E(\mathbf{Z}/\ell^{n+1})$ to $E(\mathbf{Z}/\ell^n)$, $n \geq n_0$, from which the surjectivity of the mapping from $A_{\ell^{n+1}}$ onto $A_{\ell^n}$ clearly follows, completing the proof of the lemma.

**Remark:** The above proof can be generalised to yield that for any abelian variety $A$ defined over $K$ and any point $P \in A(K)$ which does not project to a non-zero torsion point in any (geometric) subquotient of $A$, given a prime $\ell$ there are a positive density of places $v$ of $K$ such that $P$ mod $v$ has order prime to $\ell$.

**Corollary 1** $\phi(B)$ *is torsion-free where $B$ is any torsion-free subgroup of $A(K)$.*

**Proof:** Immediate from Lemma 5.

Lemma 5 allows us to strengthen Proposition 1.

**Corollary 2** *If $P$ and $Q$ are non-torsion points in $A(K)$ such that for almost all $v$ the order of $Q$ mod $v$ (i.e., the order of its image in $A_v(k_v)$) divides the order of $P$ mod $v$ (i.e., the order of its image in $A_v(k_v)$), then $Q$ is a multiple of $P$ in $A(K)$.*

**Proof:** From the conclusion of proof of Proposition 1 and the discussion after it, we get that $Q = mP + R$ for some integer $m$ and for $R \in A(K)_{tors}$. If $R$ is non-zero let the order of $R$ be divisible by a prime $\ell$. Again using our hypothesis that the order of $Q$ divides the order of $P$ in $A_v(k_v)$ for almost all $v$, we conclude that the order of $P$ in $A_v(k_v)$ is divisible by $\ell$ for almost all $v$. By Lemma 5 this can't happen and thus $Q = mP$.

Now we are in a position to complete the proof of part 2 of Theorem 1. Choose a torsion-free subgroup $B$ of $A(K)$ such that $A(K) = A(K)_{tors} \oplus B$. Then for each $P \in B$ we see from Corollary 1 that $\phi(P)$ is not torsion, and by Corollary 2 that $\phi(P) = mP$ for some $m$ that might a priori depend on $P$. Considering $P + Q$ for $P, Q \in B$ it follows that $m$ is independent of $P$. Now considering $\phi(P + P') = n(P + P')$ (for some $n$) and $\phi(P) = mP$ where $P'$ is $\in A(K)_{tors}$ and $P \in B$, we conclude that $m = n$ and thus we are done with the proof of Theorem 1.

**Remarks:**
   1. The analogs of Siegel's theorems proved by Faltings in the context of abelian varieties (i.e., Lang's conjectures on finiteness of integral points on affine subvarieties of abelian varieties that are complements of ample divisors) do not seem to yield Lemma 3 and 5 unlike in the case of elliptic curves treated in [RS].

2. Although Question 1 is posed for homomorphisms of $A(K)$, as we have seen in the proof, the heart of the question is "pointwise": namely given points $P, Q \in A(K)$, such that for almost all places $v$ of $K$ the order of $Q$ in $A_v(k_v)$ divides the order of $P$ in $A_v(k_v)$, then $Q$ is related to $P$ by an endomorphism of $A$. (Because of Lemma 6 in the next section we see that this is the essential content of the question in the case of simple abelian varities.)

## 2.2   CM elliptic curves

We begin with the following general lemma.

**Lemma 6** *Let $A$ be a finitely generated free abelian group. Let $\mathcal{O} \subset \mathcal{D}$ be an order in a division algebra $\mathcal{D}$ which is finite dimensional over $\mathbf{Q}$. Suppose $\mathcal{O}$ acts on $A$ on the left making it into a left $\mathcal{O}$-module. Suppose that $f$ is an endomorphism of $A$ as an additive group such that for all $a \in A$, there exists $f_a \in \mathcal{O}$ such that $f(a) = f_a \cdot a$. Then $f$ is multiplication by an element of $\mathcal{O}$.*

**Proof :** Clearly $A \otimes_{\mathbf{Z}} \mathbf{Q}$ is a vector space over $\mathbf{Q}$ on which $\mathcal{O} \otimes_{\mathbf{Z}} \mathbf{Q} = \mathcal{D}$ acts, making it into a $\mathcal{D}$-vector space. From the hypothesis that $f(a) = f_a \cdot a$, we find that any $\mathcal{D}$-subspace of $A \otimes_{\mathbf{Z}} \mathbf{Q}$ is stable under $f$ (extended to $A \otimes_{\mathbf{Z}} \mathbf{Q}$). Write $A \otimes_{\mathbf{Z}} \mathbf{Q} = L_1 \oplus \cdots \oplus L_n$, as a direct sum of $\mathcal{D}$-subspaces $L_i$ of dimension 1 which as has been noted is invariant under $f$, i.e., $f(L_i) \subset L_i$. Write $M_i = L_i \cap A$; then $M_i$ is a lattice in $L_i$, and $\oplus M_i$ is a subgroup of finite index in $A$. Since both $L_i$ and $A$ are invariant under $f$, so is $M_i$. Also, each $L_i$ and $A$ being invariant under $\mathcal{O}$, so is $M_i$. We will now prove that the restriction of $f$ to $M_i$ is given by multiplication by an element $f_i$ in $\mathcal{O}$.

We can clearly assume that $M_i$ is a lattice in $\mathcal{D}$ which is invariant under $\mathcal{O}$. Also, after scaling by an element of $\mathcal{D}^*$ on the right, we can assume that the lattice $M_i$ in $\mathcal{D}$ contains 1.

Suppose that $f_i(1) = \alpha_i \in \mathcal{O}$. We can thus after replacing $f_i$ by $f_i - \alpha_i$, assume that $f_i(1) = 0$. We would like to prove that $f_i$ is identically zero. Assuming the contrary, let $x$ be an element in $M_i \cap \mathcal{O}$ such that $f_i(x) \neq 0$. After scaling $x$, we can moreover assume that $f_i(x)$ belongs to $\mathcal{O}$.

By hypothesis, for every element $m \in \mathbf{Z}$, there exists an element $\lambda_m \in \mathcal{O}$ such that $f_i(m + x) = \lambda_m \cdot (m + x)$. For an element $z \in \mathcal{D}$, let $\text{Norm}(z)$ denote the determinant of the left multiplication by $z$ on $\mathcal{D}$. Since $f_i(x) = f_i(m+x) = \lambda_m \cdot (x+m)$, it follows that $\text{Norm}(m+x)$ and $\text{Norm}(f_i(x))$ are integers in $\mathbf{Q}$, and $\text{Norm}(m+x)$ divides $\text{Norm}(f_i(x))$. Since $\text{Norm}(m+x)$ is a polynomial in $m$ with coefficients in $\mathbf{Z}$ of degree equal to the dimension of $\mathcal{D}$ over $\mathbf{Q}$ of leading term 1 and constant term $\text{Norm}(x)$, the polynomial $\text{Norm}(m+x)$ as $m$ varies takes arbitrary large values, hence cannot divide the fixed integer $\text{Norm}(f_i(x))$.

We have thus proved that $f$ restricted to any 1 dimensional $\mathcal{D}$ submodule of $A \otimes_{\mathbf{Z}} \mathbf{Q}$ is multiplication by an element of $\mathcal{D}$. From this it is trivial to see that the action of $f$ on $A \otimes_{\mathbf{Z}} \mathbf{Q}$ is multiplication by an element of $\mathcal{D}$, which must moreover lie in $\mathcal{O}$, completing the proof of the lemma.

**Remark :** The lemma above holds good only in the integral version stated above, and not for vector spaces, and hence is not totally trivial. We point out an example to illustrate that the analogue of the lemma is not true for vector spaces. For this, let $K$ be a finite extension of $\mathbf{Q}$. There is an action of $K^*$ on $K$ via left or right multiplication. Let $f \in \text{Aut}_{\mathbf{Q}}(K)$ which does not arise from the action of an element of $K^*$. Such an automorphism $f$ satisfies the hypothesis of the previous lemma as for any $a \neq 0$, $f(a) \in K^*$, hence $f(a) = f_a \cdot a$ with $f_a \in K^*$, but such an $f$ does not satisfy the conclusion of the lemma.

**Corollary 3** *Question 1 has a positive answer for elliptic curves.*

Proof: This follows from the above Lemma (that we need invoke only in the case of CM elliptic curves) and Theorem 2 of [RS]. In the case of CM elliptic curves to apply the Lemma we just consider $tE(K)$ for any integer $t$ that annihilates $E(K)_{tors}$ (note that $tE(K)$ has an action of $\text{End}_K(E)$).

## 3 Morphisms between Mordell-Weil groups of two elliptic curves

As we remarked in the introduction the purely algebraic, Kummer-theoretic methods of Section 2 break down when the abelian variety $A$ is not simple and seem inadequate to answer Question 2. We discuss a approach to the question. An affirmative answer to Question 2 (in the interesting case when $\phi$ has infinite image) implies that the elliptic curves $E_i$ are isogenous over $K$ (and forces $\phi$ to have finite kernel and cokernel). Thus the question might be viewed as asking for a *Tate-like* conjecture, i.e., asking for a criterion for elliptic curves to be isogenous. The aim of the discussion is to say why it is plausible to expect an positive answer (as at the moment we are unable to do better than this!). We start off with a simple general observation.

**Proposition 2** *Let $A_1, A_2$ be abelian varieties over $K$ of dimension $g$. Then for almost all places $v$ of $K$ (more precisely, for all places $v$ of good reduction of $A$ such that $|\text{Nm}(v)| > 9.2^{4g}$):*

1. *Any injective homomorphism $\phi_v : A_1(k_v) \to A_2(k_v)$ is an isomorphism.*

2. *If $H_v$ is a subgroup of $A_1(k_v)$ and there is a surjective homomorphism $\phi_v : H_v \to A_2(k_v)$, then $H_v = A_1(k_v)$ and $\phi_v$ is an isomorphism.*

**Proof:** The proof of both statements follows from $\phi_v$ being a group homomorphism and the fact that as $A_i$'s have the same dimension, $A_i(k_v)$ have roughly the same size, i.e., the standard estimate (Weil bounds and Lefschetz fixed point formula)

$$|\text{Nm}(v)|^g - 2^{2g}|\text{Nm}(v)|^{g-1/2} \leq |A_{i,v}(k_v)| \leq |\text{Nm}(v)|^g + 2^{2g}|\text{Nm}(v)|^{g-1/2}.$$

**Discussion of Question 2:** Assume for simplicity that the $E_i$'s do not have CM and $E_2(K)$ does not have torsion (or that all of its torsion is contained in the image of $\phi$). By the arguments of [RS] we see that the image of a non-torsion point in $E_1(K)$ under $\phi$ is either 0 or again a non-torsion point. To aviod trivialities, we assume that the image of $\phi$ is infinite.

Now we make the assumption that the specialisation of $\phi(E_1(K))$ mod $v$ is all of $E_2(k_v)$ for a positive density of places $v$. Such results (i.e., elliptic analogs of Artin's conjecture on primitive roots) are conjectured by Lang and Trotter to hold in a number of situations (see [GM1] and [GM2] where some results towards this are proved).

Now from Proposition 2 and the assumption we see that the Frobenius traces $a_p(E_i)$ coincide for a positive density of primes $p$. Then by the main theorem of [Ra], as we are assuming that the $E_i$ do not have CM, the $\ell$-adic representations attached to $E_i$ are twists of each other. It is easy to see that the $\ell$-adic representations of $E_1$ and $E_2$ differ by a character that is "independent of $\ell$" (as $a_p(E_i) = 0$ for a density 0 set of primes). Thus after base change to a number field $K'$, $K'(E_1[\ell]) = K'(E_2[\ell]) = K'(E_1 \times E_2[\ell])$ for almost all $\ell$. Consider a point of infinite order $P \in E_1(K)$ such that $\phi(P) \in E_2(K)$ is non-zero (and hence of infinite order) and consider the points $(P, 0), (0, \phi(P)) \in E_1 \times E_2$. By the methods of Section 2 it again follows that $K'(\ell^{-1}(P, 0), E_1 \times E_2[\ell]) = K'(\ell^{-1}(0, \phi(P)), E_1 \times E_2[\ell])$ and using Theorem 3.2 of [JR] we get that $P$ and $\phi(P)$ are related by an isogeny from $E_1$ to $E_2$.

**Remark:** The philosophical difference which accounts for the relative ease in answering Question 1 for simple abelian varieties is in that case when proving certain "1-motives are isogenous" (which is essentially the content of Proposition 1) we assume that their abelian variety quotients are isogenous, but in Question 2 we want to prove that certain "1-motives are isogenous" without knowing *a priori* that their abelian variety quotients are isogenous.

## 4 Rigidity for arithmetic groups

We begin with the following theorem for tori.

**Proposition 3** *Given homomorphism $\phi : \mathcal{O}_K^* \to \mathcal{O}_K^*$ that reduces mod $v$ for almost all places $v$ of a number field $K$, then $\phi$ is induced by the mth power map for some integer $m$.*

**Proof:** The proof is a direct consequence of Theorem 1 of [RS] and the fact that the torsion subgroup of $K^*$ is cyclic.

We next have the following theorem for arithmetic groups.

**Theorem 2** *Let $\Gamma$ be a subgroup of $SL(2, \mathbf{Z})$ of finite index. Let $\phi$ be a non-trivial homomorphism of $\Gamma$ into itself. Assume that for all primes $p$ in an infinite set $S$ of primes, $\phi$ factors to give a homomorphism $\phi_p : SL(2, \mathbf{Z}/p) \to SL(2, \mathbf{Z}/p)$*

$$\begin{array}{ccc} \Gamma & \xrightarrow{\phi} & \Gamma \\ \downarrow & & \downarrow \\ SL(2, \mathbf{Z}/p) & \xrightarrow{\phi_p} & SL(2, \mathbf{Z}/p). \end{array}$$

*Then $\phi$ is given by the inner-conjugation action of an element in $GL(2, \mathbf{Q})$.*

**Proof:** Let $A$ be the ring which is the direct product of $\mathbf{Z}/p$ for all $p$ in $S$. Clearly $\mathbf{Z}$ is a subring of $A$, and there is thus an injective homomorphism from $SL(2, \mathbf{Z})$ to $SL(2, A)$. Since there is an injective homomorphism from $SL(2, \mathbf{Z})$ to $SL(2, A)$ for $A$, the direct product of *any* infinite set of primes, it is clear that $\phi_p$ can be trivial for at most finitely many $p$ in $S$. After replacing $S$ by this slightly smaller set, we assume that $\phi_p$ is surjective for all $p$ in $S$, and hence the $\phi_p$ are given by the inner-conjugation action of an element $g_p$ in $GL(2, \mathbf{Z}/p)$. Here we are using the well-known facts:

1. any surjective homomorphism of $SL(2, \mathbf{Z}/p)$ into itself is given by the inner-conjugation of an element of $GL(2, \mathbf{Z}/p)$.

2. any homomorphism of $SL(2, \mathbf{Z}/p)$ into itself is either trivial or is surjective if $p > 3$.

Clearly the inner-conjugation action of $g = \prod_{p \in S}(g_p) \in \prod_{p \in S} GL(2, \mathbf{Z}/p) = GL(2, A)$ takes the subgroup $\Gamma$ of $SL(2, A)$ into itself, and is the homomorphism $\phi$ when restricted to $\Gamma$.

We claim that $g$ belongs to $A^* GL(2, \mathbf{Q})$, which will complete the proof of the theorem. Observe that the inner conjugation by $g$ preserves $\Gamma$ which is a subgroup of finite index of $SL(2, \mathbf{Z})$, and hence the $\mathbf{Z}$ span of the elements of this arithmetic group which is a subring, say $R$, of $M(2, \mathbf{Z})$. It is clear that $R$ is an order inside $M(2, \mathbf{Q})$, i.e. $R \otimes \mathbf{Q} = M(2, \mathbf{Q})$. Inner conjugation by $g$ preserves this subring $R$, hence $R \otimes \mathbf{Q}$, and hence acts as an automorphism on $M(2, \mathbf{Q})$. The only automorphism of $M(2, \mathbf{Q})$ being $GL(2, \mathbf{Q})$, the homomorphism $\phi$ corresponds to an element of $GL(2, \mathbf{Q})$.

**Remarks:** 1. The proof above uses that the automorphism $\phi$ of the arithmetic group is coming from inner conjugation of an element in the larger group $GL(2, A)$ so that it preserves both additive and multiplicative structures.

2. For $SL_2$ to deduce rigidity results one just needs that the abstract homomorphism specialises for any infinite set of primes rather than for almost all or even a positive density of primes which is crucial for abelian varieties. Further unlike abelian varieties the analog of Question 1 for $SL_2$ dealt with here is not a "pointwise question" (see remark at end of Section 2.1).

3. The proof works for $SL(n, \mathbf{Z})$ for any $n$ to say that if $\phi$ is a homomorphism of a subgroup of finite index of $SL(n, \mathbf{Z})$ onto another subgroup of finite index $SL(n, \mathbf{Z})$ which specialises for infinitely many primes $p$ to give a homomorphism of $SL(n, \mathbf{Z}/p)$ to itself (we recall that an automorphism of $SL(n, \mathbf{Z}/p)$ is generated by inner automorphism from $GL(n, \mathbf{Z}/p)$, and the automorphism $A \to {}^t A^{-1}$), then $\phi$ is algebraic.

4. Because of the strong rigidity theorem, the algebraicity of abstract homomorphisms of arithmetic lattices in semi-simple Lie groups is of interest only for arithmetic lattices in $SL(2, \mathbf{R})$ which are constructed using division algebras over totally real number fields. Our method clearly applies for such lattices too.

## 5 Rational Varieties

After having considered the case of Abelian varieties and semi-simple algebraic groups, it is tempting to consider arbitrary maps of algebraic varieties over number fields which reduce nicely under reduction modulo all finite primes, i.e. such that the following diagram

$$\begin{array}{ccc} V(K) & \xrightarrow{\phi} & V(K) \\ {\scriptstyle sp_v} \downarrow & & \downarrow {\scriptstyle sp_v} \\ V(k_v) & \xrightarrow{\phi_v} & V(k_v) \end{array}$$

commutes where $k_v$ is a residue field of the ring of integers of $K$, and to ask whether such maps come from an algebraic one on $V$. The question will naturally be more meaningful if $V$ is assured of many rational points, if for instance, $V$ is a smooth projective rational variety. It seems specially interesting to investigate it for flag variety $G/P$ for a parabolic $P$ in a semi-simple split group $G$ over $K$. Here, we merely point out that the analogous question for affine line over $\mathbf{Z}$ is false, i.e. there exists a set-theoretic map from $\mathbf{Z}$ to $\mathbf{Z}$ which is not polynomial but which makes the following diagram commute.

$$\begin{array}{ccc} A^1(\mathbf{Z}) & \xrightarrow{\phi} & A^1(\mathbf{Z}) \\ {\scriptstyle sp_p} \downarrow & & \downarrow {\scriptstyle sp_p} \\ A^1(\mathbf{Z}/p) & \xrightarrow{\phi_p} & A^1(\mathbf{Z}/p). \end{array}$$

This is constructed for instance using

$$\psi(n) = a_0 + a_1 n + a_2 n(n-1) + a_3 n(n-1)(n-2) + \cdots,$$

an infinite sum, which reduces to a finite sum for each $n$, where $a_i$ are integral, and $n \geq 0$, and defining $\phi(n) = \psi(n^2)$. We find that since $\phi(n)$ is congruent to $\phi(m)$ modulo any integer $N$ for which $m$ is congruent to $n$ modulo $N$, the diagram above commutes for $\phi$.

**Remark on [BGK]:** After a preliminary form of this paper was written Ken Ribet made us aware of the preprint [BGK] which the authors then kindly made available to us. The results there have a certain amount of overlap with the proof of Theorem 1 (Proposition 2 is contained in [BGK]). They prove (Theorem E of loc. cit.) results for a larger class of abelian varieties. Using their methods and (analogs of) Lemmas 3, 5, Lemma 6 above, Question 1 can be answered for the larger class of abelian varieties covered in Thorem E of loc. cit. But the case of abelian varieties that are not simple still remains open.

*Addresses of the authors:*

CK: School of Mathematics, Tata Institute of Fundamental Research, Colaba, Bombay-400005, INDIA. e-mail: shekhar@math.tifr.res.in
*Current address:* Dept. of Math, University of Utah, 155 S 1400 E, Salt Lake City, UT 84112, USA: shekhar@math.utah.edu

DP: Harish-Chandra Research Institute, Chhatnag Road, Jhusi, Allahabad-211019, INDIA. e-mail: dprasad@mri.ernet.in